\documentclass[10pt,journal]{IEEEtran}
\usepackage{amsmath}
\usepackage{mathtools}
\usepackage{amsfonts}
\usepackage{comment}
\usepackage{subcaption}
\usepackage{amssymb}
\usepackage{bbm}
\usepackage[shortlabels]{enumitem}
\usepackage{cite}
\usepackage{colortbl}
\usepackage{tikz}
\usepackage{graphicx}
\usepackage{array}
\usepackage{multirow}
\usepackage{hhline}

\DeclareMathOperator{\col}{col}

\DeclareMathOperator{\bdiag}{diag}

\newcommand{\qand}{\quad\text{and}\quad}

\DeclarePairedDelimiter\lr{\lparen}{\rparen}

\let\leq\leqslant

\newcommand{\calA}{\ensuremath{\mathcal{A}}}
\newcommand{\calB}{\ensuremath{\mathcal{B}}}
\newcommand{\calC}{\ensuremath{\mathcal{C}}}

\newcommand{\calH}{\ensuremath{\mathcal{H}}}

\newcommand{\calM}{\ensuremath{\mathcal{M}}}
\newcommand{\calN}{\ensuremath{\mathcal{N}}}

\newcommand{\calP}{\ensuremath{\mathcal{P}}}

\newcommand{\calS}{\ensuremath{\mathcal{S}}}
\newcommand{\calT}{\ensuremath{\mathcal{T}}}

\newcommand{\calW}{\ensuremath{\mathcal{W}}}

\newcommand{\hatn}{\ensuremath{\hat{n}}}

\newcommand{\hatA}{\ensuremath{\hat{A}}}
\newcommand{\hatB}{\ensuremath{\hat{B}}}
\newcommand{\hatC}{\ensuremath{\hat{C}}}

\newcommand{\calbA}{\ensuremath{\bar{\calA}}}
\newcommand{\calbB}{\ensuremath{\bar{\calB}}}
\newcommand{\calbC}{\ensuremath{\bar{\calC}}}

\newcommand{\calhA}{\ensuremath{\hat{\calA}}}
\newcommand{\calhB}{\ensuremath{\hat{\calB}}}
\newcommand{\calhC}{\ensuremath{\hat{\calC}}}

\newcommand{\calhbA}{\ensuremath{\hat{\calbA}}}
\newcommand{\calhbB}{\ensuremath{\hat{\calbB}}}
\newcommand{\calhbC}{\ensuremath{\hat{\calbC}}}

\newcommand{\bbR}{\ensuremath{\mathbb{R}}}

\newcommand{\bmat}{\begin{matrix}}
\newcommand{\emat}{\end{matrix}}
\newcommand{\bsm}{\begin{smallmatrix}}
\newcommand{\esm}{\end{smallmatrix}}
\newcommand{\bbm}{\begin{bmatrix}}
\newcommand{\ebm}{\end{bmatrix}}
\newcommand{\bbma}{\begin{bmatrix*}[r]}
\newcommand{\ebma}{\end{bmatrix*}}
\newcommand{\bpm}{\begin{pmatrix}}
\newcommand{\epm}{\end{pmatrix}}
\newcommand{\bvm}{\begin{vmatrix}}
\newcommand{\evm}{\end{vmatrix}}
\newcommand{\bse}{\begin{subequations}}
\newcommand{\ese}{\end{subequations}}
\newcommand{\beq}{\begin{equation}}
\newcommand{\eeq}{\end{equation}}
\newcommand{\beqn}{\begin{equation*}}
\newcommand{\eeqn}{\end{equation*}}

\newcommand{\ben}{\renewcommand{\labelenumi}{\arabic{enumi}.}
\renewcommand{\theenumi}{\arabic{enumi}}\begin{enumerate}}

\newcommand{\een}{\end{enumerate}}

\newcommand{\beni}{\renewcommand{\labelenumi}{\roman{enumi}.}
\renewcommand{\theenumi}{\roman{enumi}}\begin{enumerate}}

\newcommand{\eeni}{\end{enumerate}}

\newcommand{\bena}{\renewcommand{\labelenumi}{\alph{enumi}.}
\renewcommand{\theenumi}{\alph{enumi}}\begin{enumerate}}

\newcommand{\eena}{\end{enumerate}}
\newcommand{\bit}{\begin{itemize}}
\newcommand{\eit}{\end{itemize}}

\newtheorem{thm}{Theorem}
\newtheorem{defn}[thm]{Definition}
\newtheorem{lem}[thm]{Lemma}

\newtheorem{cor}[thm]{Corollary}
\newtheorem{example}[thm]{Example}
\newtheorem{remark}[thm]{Remark}
\newtheorem{prb}[thm]{Problem}
\newtheorem{prp}[thm]{Proposition}

\newcommand{\bthm}{\begin{thm}}
\newcommand{\ethm}{\end{thm}}
\newcommand{\blem}{\begin{lem}}
\newcommand{\elem}{\end{lem}}
\newcommand{\bprop}{\begin{prp}}
\newcommand{\eprop}{\end{prp}}
\newcommand{\bex}{\begin{example}}
\newcommand{\eex}{\end{example}}
\newcommand{\bas}{\begin{assumption}}
\newcommand{\eas}{\end{assumption}}
\newcommand{\bre}{\begin{remark}}
\newcommand{\ere}{\end{remark}}
\newcommand{\bcor}{\begin{cor}}
\newcommand{\ecor}{\end{cor}}
\newcommand{\bdfn}{\begin{defn}}
\newcommand{\edfn}{\end{defn}}
\newcommand{\bcon}{\begin{conjecture}}
\newcommand{\econ}{\end{conjecture}}
\newcommand{\bali}{\begin{aligned}}
\newcommand{\eali}{\end{aligned}}
\newcommand{\eprb}{\end{prb}}
\newcommand{\bprb}{\begin{prb}}

\newcommand{\nset}[1]{\ensuremath{\{1,\ldots,#1\}}}

\newcommand{\Ci}{\textup{(Ci)}}
\newtheorem{corollary}[thm]{Corollary}
\newtheorem{proposition}[thm]{Proposition}
\newtheorem{problem}[thm]{Problem}
\newtheorem{assumption}[thm]{Assumption}

\usepackage{tcolorbox}
\usetikzlibrary{shapes.geometric, arrows, automata}
\usepackage[colorlinks=true,linkcolor=green]{hyperref}
\hypersetup{
	colorlinks,
	linkcolor={blue!80!black},
	citecolor={green!80!black},
	urlcolor={brown!80!black}
}
\def\equationautorefname~#1\null{Equation~(#1)\null}

\newcommand{\magenta}{\textcolor{magenta}}

\usetikzlibrary{calc,patterns,decorations.pathmorphing,decorations.markings}

\newcounter{todocounter}
\usepackage[colorinlistoftodos]{todonotes}

\makeatletter

\newcommand{\Rmnum}[1]{\expandafter\@slowromancap\romannumeral #1@}
\makeatother

\usepackage[normalem]{ulem} 

\usepackage{pmat}

\begin{document}
%
\title{Scalable controllability analysis of structured networks}


\author{\IEEEauthorblockN{J. Jia,
B. M. Shali,
H. J. van Waarde,
M. K. Camlibel,~\IEEEmembership{Member,~IEEE},
and
H. L. Trentelman,~\IEEEmembership{Fellow,~IEEE}
}%
\thanks{
    J. Jia is with the Key Laboratory of Advanced Process Control for Light Industry (Ministry of Education), Jiangnan University, Wuxi 214122, China; Email: \texttt{jiajiajia0218@163.com}.

    B. M. Shali, H. J. van Waarde, M. K. Camlibel and H. L. Trentelman are with the Bernoulli Institute for Mathematics, Computer Science and Artificial Intelligence, University of Groningen, The Netherlands; Email: {\tt b.m.shali@rug.nl}, {\tt h.j.van.waarde@rug.nl},	{\tt m.k.camlibel@rug.nl}, {\tt h.l.trentelman@rug.nl}.

    This work is partially supported by Jiangsu Provincial Natural Science Foundation of China (BK20201340), China Postdoctoral Science Foundation (2018M642160) and the 111 Project (B12018).
    }}
%

%



\IEEEtitleabstractindextext{%
%
%
%
%
%

\begin{abstract}
This paper deals with strong structural controllability of structured networks.
A structured network is a family of structured systems (called node systems) that are interconnected by means of a structured interconnection law.
The node systems and their structured interconnection law are given by pattern matrices.
It is shown that a structured network is strongly structurally controllable if and only if an associated structured system is.
This structured system will in general have a very large state space dimension, and therefore existing tests for verifying strong structural controllability are not tractable.
The main result of this paper circumvents this problem. We show that controllability can be tested by replacing the original network by a new network in which all original node systems have been replaced by (auxiliary) node systems with state space dimensions either 1 or 2.
Hence, controllability of the original network can be verified by testing controllability of a structured system with state space dimension at most twice the number of node systems, regardless of the state space dimensions of the original node systems.

\end{abstract}

}

\maketitle

\IEEEdisplaynontitleabstractindextext

%
\IEEEpeerreviewmaketitle

\maketitle
\section{Introduction}
This paper deals with strong structural controllability of interconnections of structured systems.
The starting point is a collection of linear structured input-state-output systems, called the {\em node systems}.
These are systems in which the system matrices are not given by matrices with real entries, but, instead, by so-called pattern matrices.
These pattern matrices indicate which entries in the system matrices are equal to zero, which are arbitrary nonzero, and which are completely arbitrary (zero or nonzero).
In addition, a structured interconnection topology is given, also in terms of pattern matrices.
Such structured interconnection topologies enable us to make a distinction between links that are certainly present, and links that might be present or not.
By formally interconnecting the node systems through their inputs and outputs as prescribed by the structured interconnection topology, and at the same time specifying a new external control input, we obtain a new (high dimensional) structured system.
This system will be called a {\em structured network}.
In this paper we will deal with finding conditions on the interplay between the node systems and interconnection topology such that this structured network is strongly structurally controllable.
This means that for all particular choices of node systems with the given structure, and all particular choices of the interconnection topology, the resulting interconnection is controllable in the classical sense.

The problem of finding conditions for controllability of interconnected systems has been studied before, mainly in the context of (non-structured) node systems and interconnection topologies represented by real matrices. Here, we refer to early work by Gilbert \cite{gilbert} dealing with controllability of systems in parallel, series and feedback interconnections, and related work by Callier and Nahum \cite{callier}. More recent references on controllability of networked system are the work of Fuhrmann and Helmke \cite{fuhrmann} and Hara et al.  \cite{hara}. We also refer to \cite{TT2018} and \cite{tanner2004, Zhou2015, WWC2017, WCWT2016, RJME2009}.
While these references all deal with interconnections of systems represented by numerical matrices, in the present paper we deal with controllability of structured interconnections of structured systems. Up to now, the major part of related work  has been done for the special case that all node systems are single integrators. Structured interconnections of such single integrators can be considered as structured systems themselves, with the system matrices given by pattern matrices. Structural controllability of such systems has been studied extensively in \cite{LSB2011} and, e.g.,\cite{TD2015,Jia2020,RA2013,MHM2018,JTBC2018,CM2013,MZC2014}. Different concepts of pattern matrices exist. Traditionally they have been defined as matrices with two kinds of entries, namely either $0$ or a (zero or nonzero) indeterminate, see \cite{L1974, DCW2003, LSB2011, C2019} and the references therein. More recently, in \cite{PPKAL2019} and \cite{Jia2020}, a more general kind of pattern matrices has been introduced, allowing three kinds of entries, namely $0$, a nonzero indeterminate, and an indeterminate that can be either zero or nonzero. It is this latter concept of pattern matrix that will be used in this paper.

A structured system is called {\em weakly} ({\em strongly}) structurally controllable if for almost all (for all) possible choices of values of the indeterminate parameters, the corresponding (numerical) linear systems are controllable.
Conditions for weak and strong structural controllability have been provided entirely in terms of the graph associated with the structured interconnection, using concepts like cactus graphs\cite{L1974}, maximal matchings \cite{LSB2011}, constrained matchings \cite{CM2013}, zero forcing sets \cite{MZC2014,TD2015} and color change rules \cite{Jia2020}.

As outlined earlier in this introduction, in the present paper we will study structural controllability of structured networks in which the node systems are structured systems themselves (in contrast to single integrators).
Related work, albeit on {\em weak} structural controllability of structured networks, can be found in e.g., \cite{CCVFT2012, WWC2017, CK2019}.
In this paper, instead, we will study {\em strong} structural controllability of structured networks.

Since networks may in general consist of a large number of node systems, with each of these having a  possibly large state space dimension themselves, the global network dimension may become prohibitively large.
Therefore, it is important to establish so-called {\em scalable} methods to determine whether a structured network is controllable. In the present paper we will indeed develop a test for controllability of structured networks in which the complexity is independent of the state space dimensions of the node systems.

The main contributions of this paper are the following.
\begin{enumerate}
\item
We show that a given structured network is (strongly structurally) controllable if and only if an associated structured system is. In principle, this result makes it possible to apply existing tests as cited above to check controllability of structured networks.
\item
In order to reduce the complexity of the previous tests, we show that controllability of a given structured network can be tested by replacing the original network by a new network in which all node systems have been replaced by (auxiliary) standard node systems with state space dimension either 1 or 2.
This means that controllability of any network can be checked by testing controllability of a structured system of state space dimension at most twice the number of node systems, regardless of the dimensions of the original node systems.
\item
In order to perform the above reduction step, we introduce a graph theoretic test (in terms of a color change procedure) to determine by which auxiliary first or second order standard node system any of the original node systems should be replaced.
\end{enumerate}

The outline of this paper is as follows.
In Section \ref{S:pre}, we introduce some notation, and review basic material on pattern matrices. In Section \ref{sec:NSS}, we explain what we mean by a structured network and formulate the main problem that will be considered in this paper. We also present an example of a structured network, which will be a running example throughout the paper.
In Section~\ref{s:condition_for_controllability}, we first review addition and multiplication of pattern matrices.
Next, we show that controllability of a structured network with node systems from a particular class of single-input single-output systems can be rephrased as controllability of an ordinary structured system.
Section~\ref{s:reduce_algebraic} is the key section of this paper.
It is shown that controllability of a structured network is equivalent to controllability of a new network in which all original node systems have been replaced by new node systems from a finite set of standard node systems with state space dimension either 1 or 2.
In Section~\ref{s:reduce_graph}, we establish a graph theoretic tool for determining which one of the standard node systems should replace a given original node system.
Finally,  in Section~\ref{s:s5_con}, we formulate our conclusions.
All proofs in this paper have been  defered to the Appendix.

\section{Preliminaries} \label{S:pre}

Given a set of matrices $\{ A_1, \ldots, A_n\}$, we denote
\[
\bdiag(A_1, \ldots, A_n) = \bbm A_1 & \cdots & 0\\ \vdots & \ddots & \vdots \\ 0& \cdots & A_n\ebm.
\]
In addition, if these matrices have the same column dimensions, we denote
\[
\col(A_1, \ldots, A_n) = \bbm A_1 \\ \vdots\\ A_n\ebm.
\]

An important role in this paper is played by {\em pattern matrices}.
These are matrices with entries in the set of symbols $\{0, \ast, ?\}$.
The set of all $p \times q$ pattern matrices is denoted by $\{0,\ast,?\}^{p \times q}$. For a given pattern matrix $\calM \in \{0,\ast,?\}^{p \times q}$, we define the {\em pattern class} of $\calM$ as the subset of $\bbR^{p \times q}$ given by
\[
\bali   
\calP(\calM) = \{M \in \bbR^{p \times q} \mid M_{ij} = 0 \mbox{ if } \calM_{ij}= 0,&\\ M_{ij} \neq 0 \mbox{ if } \calM_{ij} = \ast&\}.
\eali
\]
This means that for a given matrix $M \in \calP(\calM)$, the entry $M_{ij}$ has the real value $0$ if $\calM_{ij} = 0$, is a nonzero real number if $\calM_{ij} = \ast$, and is an arbitrary real number if $\calM_{ij} = ?$.
A pattern matrix $\calM \in \{0,\ast,?\}^{p \times q}$ with $p \leq q$ is said to have {\em full row rank} if $M$ has full row rank for every $M \in \calP(\calM)$.

\section{Structured Networks} \label{sec:NSS}
We will now first review the concept of structured system. Subsequently, we will define what we mean by a structured network. For given pattern matrices $\calA$, $\calB$ and $\calC$ of dimensions $n \times n$, $n \times m$ and $p \times n$ respectively, we define the \emph{structured system} associated with these pattern matrices as the family of linear time-invariant systems
\begin{align}
\dot{x} &= A x + Bu, \label{issys} \\
y & = C x,
\end{align}
where $A \in \calP(\calA)$, $B \in \calP(\calB)$ and $C \in \calP(\calC)$. This structured system will be denoted by $(\calA, \calB, \calC)$. Similarly, the family of systems \eqref{issys} is denoted by $(\calA, \calB)$. We will say that $(\calA, \calB)$ is \emph{strongly structurally controllable} (or simply, controllable) if \eqref{issys} is controllable for all  $A \in \calP(\calA)$ and $B \in \calP(\calB)$.

In this paper, we will study {\em structured networks}. Such a network is a family of structured systems that are interconnected by means of a structured interconnection law. More specifically, assume that for $k = 1,2, \ldots,N$ we have structured systems $(\calA_k, \calB_k,\calC_k)$,
where $\calA_k$ has dimensions $n_k \times n_k$, $\calB_k$ is $n_k \times r_k$, and $\calC_k$ is $p_k \times n_k$. This will be called the {\em structured node system} at node $k$. Define $r := \sum_{k = 1}^N r_k$ and $p := \sum_{k = 1}^N p_k$. Next, a {\em structured interconnection law} is given by an $r \times p$ block pattern matrix
\begin{equation}
\label{eqW}
\calW = \bbm \calW_{11} &  \ldots & \calW_{1N} \\
               \vdots & \ddots  & \vdots \\
             \calW_{N1} &  \ldots & \calW_{NN} \ebm
\end{equation}
and an $r \times m$ block pattern matrix
\begin{equation}
\label{eqH}
\calH = \bbm \calH_1 \\ \vdots \\ \calH_{N} \ebm.
\end{equation}
 The corresponding structured network is now defined as the family of networks obtained by interconnecting $N$ node systems
\begin{equation} \label{e:agentk}
\bali
\dot{x}_k &= A_k x_k + B_kv_k, \\
y_k & = C_k x_k,
\eali
\end{equation}
with $(A_k,B_k,C_k) \in \mathcal{P}(\calA_k) \times \calP(\calB_k) \times \calP(\calC_k)$,
using an interconnection law
\begin{equation} \label{e:law}
v_k = \sum_{j=1}^{N} W_{kj} y_j +  H_k u,
\end{equation}
with $W_{kj} \in \calP(\calW_{kj})$ and $H_k \in \calP(\calH_{k})$. The new variable $u$ is an external control input taking its values in $\mathbb{R}^m$. By introducing the block diagonal matrices
\beq \label{eq:ABC}
\bali
A &= \bdiag(A_1,\ldots,A_N),\\
B &= \bdiag(B_1,\ldots,B_N),\\
C &= \bdiag(C_1,\ldots,C_N),\\
\eali
\eeq
the interconnection of \eqref{e:agentk} and \eqref{e:law} can be represented compactly as
\beq \label{eq:Snetworksystem}
\dot{x}  = (A+ B W C) x + B H u.
\eeq
Here, $x = \col(x_1, \ldots, x_N)$ denotes the vector obtained by stacking the states of all node systems. Obviously, $x  \in \mathbb{R}^n$ with $n := \sum_{k=1}^N n_k$.

Now introduce the block pattern matrices
\begin{equation}
\label{eq:blkdiagpatterns}
\begin{aligned}
    \calA &= \bdiag(\calA_1, \ldots, \calA_N), \\
	\calB &= \bdiag(\calB_1, \ldots, \calB_N), \\
	\calC &= \bdiag(\calC_1, \ldots, \calC_N).
\end{aligned}
\end{equation}
It is then clear that our structured network consists of all systems  \eqref{eq:Snetworksystem}, where $A \in \calP(\calA)$, $B \in \calP(\calB)$, $C \in \calP(\calC)$, $W \in \calP(\calW)$ and $H \in \calP(\calH)$.
This structured network will be denoted by $(\calA,\calB,\calC,\calW,\calH)$.

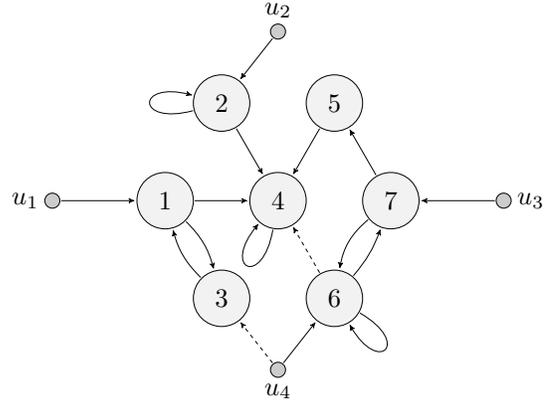
\begin{figure}[h!]
		\centering
        \scalebox{0.5}{
            \begin{tikzpicture}[->,>=stealth',shorten >=1pt,auto,node distance=2.8cm, semithick]
                \newcommand{\sR}{.8cm}
                \newcommand{\bR}{3cm}

                \tikzset{VS1/.style = {shape = circle,
                        color=black,
                        fill=white!80!black,
                        minimum size=0.4cm,
                        text = black,
                        inner sep = 2pt,
                        outer sep = 1pt,
                        draw}
                }

                \tikzset{VS2/.style = {shape = circle,
                        color=black,
                        fill=white!95!black,
                        minimum size=1.5cm,
                        text = black,
                        inner sep = 2pt,
                        outer sep = 1pt,
                        draw}
                }

                \node[VS2] (1) at (180: \bR) {\scalebox{2}{$1$}};
                \node[VS2] (2) at (120: \bR) {\scalebox{2}{$2$}};
                \node[VS2] (3) at (240: \bR) {\scalebox{2}{$3$}};
                \node[VS2] (4) at (0,0) {\scalebox{2}{$4$}};
                \node[VS2] (5) at (60: \bR) {\scalebox{2}{$5$}};
                \node[VS2] (6) at (300: \bR) {\scalebox{2}{$6$}};
                \node[VS2] (7) at (0: \bR) {\scalebox{2}{$7$}};

                \node[VS1] (u1) at (180: 2*\bR) [label=left:\scalebox{2}{$u_1$}] {};
                \node[VS1] (u2) at (90: 1.5*\bR) [label=above:\scalebox{2}{$u_2$}] {};
                \node[VS1] (u3) at (0: 2*\bR) [label=right:\scalebox{2}{$u_3$}] {};
                \node[VS1] (u4) at (-90: 1.5*\bR) [label=below:\scalebox{2}{$u_4$}] {};


                \draw (1) to [out=-45, in=105, looseness = 1] (3);
                \draw (1) to (4);
                \draw (2) to (4);
                \draw (3) to [out=135, in=-75, looseness = 1] (1);
                \draw (5) to (4);
                \draw[dashed] (6) to (4);
                \draw (6) to [out=50, in=-110, looseness = 1] (7);
                \draw (7) to (5);
                \draw (7) to [out=-140, in=80, looseness = 1] (6);

                \draw (u1) to (1);
                \draw (u2) to (2);
                \draw (u3) to (7);
                \draw[dashed] (u4) to (3);
                \draw (u4) to (6);

                \draw (2) to [out=-165, in=165, looseness = 10] (2);
                \draw (4) to [out=-100, in=-130, looseness = 10] (4);
                \draw (6) to [out=-30, in=-60, looseness = 10] (6);
            \end{tikzpicture}
        }
		\caption{The structured network of Example \ref{ex:1}.}
        \vspace{-4mm}
		\label{f:topology}
	\end{figure}

In this paper we are interested in controllability of this network. We will say that $(\calA,\calB,\calC,\calW,\calH)$ is \emph{strongly structurally controllable} if  \eqref{eq:Snetworksystem} is controllable for all  $A \in \calP(\calA)$, $B \in \calP(\calB)$, $C \in \calP(\calC)$, $W \in \calP(\calW)$ and $H \in \calP(\calH)$. In this paper we will then simply call the structured network controllable.
 The problem that we will investigate is the following:
 \vspace{2mm}
\begin{problem}
Find necessary and sufficient conditions under which the structured network $(\calA,\calB,\calC,\calW,\calH)$ is controllable.
\end{problem}
\vspace{2mm}
\bex \label{ex:1}
We will illustrate the set up introduced above using a network of mechanical systems. This example will be the leading example throughout this paper. Specifically, we consider a network consisting of $7$ structured single-input single-output node systems and $4$ external inputs interconnected through the structured interconnection law defined by the pattern matrices
\begin{equation} \label{eq:intstruc_W}
\calW  :=
\bbm
0       &  0        &   \ast        &  0       &  0       &   0     & 0   \\
0       &  ?        &   0         &  0      &  0       &   0     & 0    \\
\ast   &  0       &  0     &  0      &  0       &   0     & 0    \\
\ast   &  \ast  &   0          &  \ast    &  \ast   &    ?    & 0    \\
0       &  0       &   0         &  0  &  0       &  0       & \ast \\
0       &  0      &   0          &  0      &  0       & \ast    & \ast \\
0       &  0      &   0          &  0      &  0       & \ast      & 0     \\ \ebm,
\calH :=
\bbm
\ast  &  0        &   0      & 0\\
0      &  \ast        &   0    & 0  \\
0      &  0       &   0   & ?\\
0      &  0       &   0    & 0    \\
0      &  0       &   0      &0 \\
0      &  0      &   0      &\ast \\
0      &  0      &   \ast    & 0   \\
\ebm.
\end{equation}
A graphical representation of this structured interconnection law is depicted in Figure~\ref{f:topology}.

Suppose that the $k$-th node system has actuated mass-spring-damper dynamics of the form
\beq
f_k = m_k \ddot{p}_k + c_k \dot{p}_k + \ell_k p_k,
\eeq
in which $m_k$, $c_k$ and $\ell_k$ denote mass, damper constant and spring constant, respectively, and $f_k$ and $p_k$ represent force and position.
By introducing $\xi_k^\top =\bbm p_k & \dot{p}_k \ebm$, we obtain
\beq
\begin{aligned}
\dot{\xi}_k  &= \underbrace{\bbm 0 & 1\\ -\frac{\ell_k}{m_k} & -\frac{c_k}{m_k}\ebm}_{:= X_k} \xi_k + \underbrace{\bbm 0 \\ \frac{1}{m_k}\ebm}_{:= Y_k} f_k,\\
z_k  & = \underbrace{\bbm 1 & 0 \ebm}_{: =  Z_k} \xi_k, \\
\end{aligned}
\eeq
where $z_k$ is the output reflecting that position is measured.
Each node system is actuated by a dynamic controller of the form
\beq \label{eq:controller}
\begin{aligned}
\dot{\omega}_k &= K_k  \omega_k + L_k z_k,\\
f_k & = M_k \omega_k + v_k,
\end{aligned}
\eeq
with $K_k \in \bbR^{2 \times 2}$, $L_k \in \bbR^{2}$ and $M_k \in \bbR^{1 \times 2}$.
By defining $x_k := \col (\xi_k,\omega_k)$ and taking the output $y_k$ of the controlled node system equal to $z_k$, we obtain the resulting node dynamics
\beq \label{eq:node_sys}
\begin{aligned}
\dot{x}_k &= A_k  x_k + B_{k} v_k,\\
y_k & = C_{k} x_k,
\end{aligned}
\eeq
where
$$
A_k = \bbm X_k & Y_k M_k \\ L_k Z_k & K_k \ebm,\quad B_k = \bbm Y_k \\ 0 \ebm,\quad C_k = \bbm Z_k & 0 \ebm.
$$
Suppose now that the nonzero parameters $m_k, c_k$ and $\ell_k$ are not known exactly. This means that the matrices $X_k$, $Y_k$ and $Z_k$ are not known exactly. Assume that the controller for the first node system is chosen as
\[
K_1= \bbm 0 & k_1\\ 0 & 0\ebm, \quad L_1 = \bbm 0 \\ l_1\ebm \qand M_1 = \bbm m_1 & 0 \ebm,
\]
where $k_1$, $l_1$ and $m_1$ are nonzero real parameters. Then the first node system can be represented as the structured system $(\calA_1,\calB_1, \calC_1)$, with
\[
\calA_1 =
\bbm
0 & \ast & 0 & 0\\
\ast & \ast & \ast & 0\\
0& 0 &0 & \ast\\
\ast & 0 &0 & 0\\
\ebm,
\quad \calB_1 = \bbm 0 \\ \ast \\ 0 \\ 0\ebm, \quad \calC_1^\top = \bbm \ast \\ 0 \\ 0 \\ 0 \ebm.\]
Likewise, assume that controllers for the other node systems have been chosen resulting in
$$
\begin{aligned}
 \calA_1 = \calA_4 = \calA_7  &=
\bbm
0 & \ast & 0 & 0\\
\ast & \ast & \ast & 0\\
0& 0 &0 & \ast\\
\ast & 0 &0 & 0\\
\ebm, \quad
\calA_2 =
\bbm
0 & \ast & 0 & 0\\
\ast & \ast & ? & \ast \\
\ast& 0 &0 & 0\\
0 & 0 & \ast & 0\\
\ebm,\\
 \calA_3  = \calA_5&=
\bbm
0      & \ast & 0    & 0 \\
\ast & \ast & 0     & 0 \\
\ast     & 0     & 0     & \ast \\
0     & 0     & \ast & 0
\ebm, \quad
\calA_6  =
\bbm
0          & \ast & 0     & 0 \\
\ast     & \ast & 0      & 0 \\
0       & 0     & \ast     & \ast \\
\ast     & 0     & 0     & 0 \\
\ebm.
\end{aligned}
$$
The matrices $B_k$ and $C_k$ in \eqref{eq:node_sys} for $k = 2,3, \ldots, 7$ remain to have the structure $\calB_k = \calB_1$ and  $\calC_k = \calC_1.$
The entire network is now described by the 5-tuple $(\calA,\calB,\calC,\calW,\calH)$ with $\calA,\calB$ and $\calC$ defined as in \eqref{eq:blkdiagpatterns}, and $\calW$ and $\calH$ defined by \eqref{eq:intstruc_W}.
\eex

\section{Algebraic conditions for controllability of structured networks}\label{s:condition_for_controllability}

In this section, we will provide conditions under which $(\calA,\calB,\calC,\calW,\calH)$ is controllable. In view of Equation \eqref{eq:Snetworksystem}, it may be tempting to characterize controllability of a networked structured system by analyzing the controllability of some sort of structured system of the form $(\mathcal{A}+\mathcal{B}\mathcal{W}\mathcal{C},\mathcal{B}\mathcal{H})$. However, note that the matrices $\mathcal{A}+\mathcal{B}\mathcal{W}\mathcal{C}$ and $\mathcal{B}\mathcal{H}$ are composed of sums and products of pattern matrices, \emph{which have not been formally defined yet}.

Before presenting our controllability results, we will therefore first recall the notions of addition and multiplication of pattern matrices \cite{Shali2020}. First, addition and multiplication of the symbols $0,\ast$ and $?$ are defined in Table \ref{ta:results}.

\begin{table}[h!]
	\centering
	\captionsetup{width=.75\textwidth, justification = centering}
	\caption{Addition  and multiplication within the set $\{0,\ast,?\}$.}
	\begin{tabular}{c|ccc}
		$+$ & $0$ & $\ast$ & $?$  \\  \hline
		$0\rule{0pt}{2.2ex}$ & $0$ & $\ast$ & ${?}$ \\
		$\ast$& $\ast$ & $?$ & $?$  \\
		$?$& $?$ & $?$ & $?$
	\end{tabular}
	\qquad
	\begin{tabular}{c|ccc}
		$\boldsymbol{\cdot}$  & $0$ & $\ast$ & $?$ \\  \hline
		$0\rule{0pt}{2.2ex}$ &  $0$ & $0$ & $0$ \\
		$\ast$& $0$ & $\ast$ & $?$ \\
		$?$&  $0$ & $?$ & $?$
	\end{tabular}
	\label{ta:results}
\end{table}
Based on the operations defined in this table, addition of pattern matrices is  then defined as follows.
\bdfn \label{d:mPa}
Let  $\calM, \calN \in \{0,\ast,?\}^{p \times q}$.
The sum of these pattern matrices, $\calM + \calN \in \{0,\ast,?\}^{p \times q}$, is defined as
\beqn
(\calM + \calN)_{ij}:=  \calM_{ij} + \calN_{ij}.
\eeqn
\edfn
We define $\calP(\calM) + \calP(\calN)$ as the usual Minkowski sum of sets, that is,
$$
\calP(\calM) + \calP(\calN) := \{M+N \mid M \in \calP(\calM) \mbox{ and }  N \in \calP(\calN)\}.
$$
We now have the following proposition.
\bprop\label{p:1} \cite[Proposition 1]{Shali2020}
For pattern matrices $\calM$ and $\calN$ of the same dimensions,
$\calP(\calM) + \calP(\calN) = \calP(\calM + \calN)$.
\eprop

In addition, we define multiplication of pattern matrices.
\bdfn \label{d:mPm}
Let  $\calM \in \{0,\ast,?\}^{p \times q}$ and $\calN \in \{0,\ast,?\}^{q \times s}$.
Then the product $\calM \calN \in \{0,\ast,?\}^{p \times s}$ is defined by
\beqn
(\calM \calN)_{ij} :=  \sum_{\ell = 1}^{q} \calM_{i\ell} \boldsymbol{\cdot} \calN_{\ell j}.
\eeqn
\edfn

In addition, we define
$$
\calP(\calM)\calP(\calN) := \{MN \mid M \in \calP(\calM) \mbox{ and }  N \in \calP(\calN)\}.
$$
It is known that the equality $\calP(\calM)\calP(\calN) = \calP(\calM \calN)$ \emph{does not hold} for general pattern matrices $\mathcal{M}$ and $\mathcal{N}$ \cite[Example 1]{Shali2020}. Nonetheless, as we demonstrate next, such an equality can be derived if $\mathcal{M}$ and $\mathcal{N}$ have a special structure.

\blem \label{l:PEqual}
Consider two pattern matrices $\calM \in \{0,\ast,?\}^{p \times q}$ and $\calN \in \{0,\ast,?\}^{q \times s}$.
Then, the equality
$$
\calP(\calM)\calP(\calN) = \calP(\calM\calN)
$$
holds if at least one of the following two conditions holds:
\begin{enumerate}[(i)]
\item  each row of $\calN$ has exactly one entry equal to $\ast$ and the remaining entries are zero, \label{i:PEqual_1}
\item each column of $\calM$ has exactly one entry equal to $\ast$ and the remaining entries are zero. \label{i:PEqual_2}
\end{enumerate}
	\elem
	\begin{IEEEproof}
		We will only consider the case that \ref{i:PEqual_1} holds.
		The other case follows from the fact that
		$\calP(\calM\calN) = \calP(\calN^\top \calM^\top)^\top$.
		To begin with, denote by $\calM_i$ the $i$-th column of $\calM$ and by $\calN_i$ the $i$-th row of $\calN$ for $i = 1, \ldots, q$.
Then we can write
\begin{equation*}
\begin{aligned}
&\calP(\calM\calN)  = \calP\left(\sum_{i=1}^{q} \calM_i \calN_i\right) = \sum_{i=1}^{q} \calP(\calM_i\calN_i),\\
&\calP(\calM)\calP(\calN)  = \sum_{i=1}^{q} \calP(\calM_i)\calP(\calN_i).
\end{aligned}
\end{equation*}
Since $\calN_i$ has exactly one entry which is $\ast$ and the remaining entries are zero, it follows that
$\calP(\calM_i\calN_i) =  \calP(\calM_i)\calP(\calN_i).$
Thus, we have
\begin{equation*}
			\calP(\calM\calN) = \calP(\calM)\calP(\calN).
\end{equation*}
This completes the proof.
\end{IEEEproof}

We now make the following two simplifying assumptions that will be in place in the rest of the paper.

\begin{assumption}
\label{asmSISO}
The node systems are \emph{single-input single-output}, i.e., $\mathcal{B}_k \in \{0,*,?\}^{n_k \times 1}$ and $\mathcal{C}_k \in \{0,*,?\}^{1 \times n_k}$ for all $k = 1,2,\dots,N$.
\end{assumption}

\begin{assumption}
\label{asmAST}
For all $k = 1,2,\dots,N$, each entry of the vectors $\mathcal{B}_k, \mathcal{C}_k^\top \in \{0,*,?\}^{n_k \times 1}$ is zero, except for one entry that is equal to $\ast$.
\end{assumption}

By Assumptions~\ref{asmSISO} and \ref{asmAST}, the matrix $\calB$ (as defined in \eqref{eq:blkdiagpatterns}) is such that the entries of every column are equal to zero except for one entry equal to $\ast$. Similarly, $\calC$ is such that the entries of every row are zero except for one $\ast$ entry. This special structure will allow us to apply Lemma~\ref{l:PEqual}. The following theorem is the main result of this section and characterizes controllability of $(\calA,\calB,\calC,\calW,\calH)$.

\bthm \label{t:GGC}
The structured network $(\calA,\calB,\calC,\calW,\calH)$ is controllable if and only if $(\calA + \calB \calW \calC, \calB \calH )$ is controllable.
\ethm

\begin{IEEEproof}
By Proposition~\ref{p:1} we have that $$\mathcal{P}(\mathcal{A}+\mathcal{B}\mathcal{W}\mathcal{C}) = \mathcal{P}(\mathcal{A}) + \mathcal{P}(\mathcal{B}\mathcal{W}\mathcal{C}).$$
Furthermore, by the special structure of the pattern matrices $\mathcal{B}$ and $\mathcal{C}$, Lemma~\ref{l:PEqual} implies that
\begin{align*}
    \mathcal{P}(\mathcal{A}+\mathcal{B}\mathcal{W}\mathcal{C}) &= \mathcal{P}(\mathcal{A}) + \mathcal{P}(\mathcal{B})\mathcal{P}(\mathcal{W})\mathcal{P}(\mathcal{C}), \\
    \mathcal{P}(\mathcal{B}\mathcal{H}) &= \mathcal{P}(\mathcal{B})\mathcal{P}(\mathcal{H}).
\end{align*}
As such, $(\calA,\calB,\calC,\calW,\calH)$ is controllable if and only if $(\calA + \calB \calW \calC, \calB \calH )$ is controllable, which proves the theorem.
\end{IEEEproof}

Theorem~\ref{t:GGC} is relevant because it relates the controllability of a structured network to that of an ordinary (albeit large) structured system, whose controllability properties are well understood \cite{Jia2020}.
In fact, we recall the following result \cite[Theorem 6]{Jia2020} that relates controllability of an arbitrary structured system $(\calA,\calB)$ to the full row rank properties of two pattern matrices.
Before stating this result, we introduce some notation. For a given square pattern matrix $\mathcal{X} \in \{0,*,?\}^{n \times n}$, we define $\bar{\mathcal{X}}$ as the pattern matrix obtained from $\mathcal{X}$ by modifying its diagonal entries as follows:
$$
\bar{\mathcal{X}}_{ii} = \begin{cases}
* & \text{if } \mathcal{X}_{ii} = 0, \\
? & \text{otherwise}.
\end{cases}$$
\begin{proposition}\label{p:fullrank}
The system $(\calA,\calB)$ is controllable if and only if both $\begin{bmatrix}\mathcal{A} & \calB \end{bmatrix}$ and $\begin{bmatrix}\bar{\mathcal{A}} & \calB \end{bmatrix}$ have full row rank.
\end{proposition}

Using this result, we immediately obtain the following algebraic characterization of network controllability.

\begin{proposition}\label{p:algebraic}
The structured network $(\calA,\calB,\calC,\calW,\calH)$ is controllable if and only if the pattern matrices $\begin{bmatrix} \calA + \calB\calW\calC & \calB\calH \end{bmatrix}$ and $\begin{bmatrix} \bar{\calA} + \calB\calW\calC & \calB\calH \end{bmatrix}$ have full row rank.
\end{proposition}

\begin{IEEEproof}
The proof follows immediately from Theorem~\ref{t:GGC} and Proposition~\ref{p:fullrank} by noting that $\overline{\mathcal{A}+\calB \calW \calC} = \bar{\mathcal{A}}+\calB \calW \calC$.
\end{IEEEproof}

Note that in the special case that the node dynamics are ``single integrators", i.e., if $\mathcal{A}_k = 0$ and $\mathcal{B}_k = \mathcal{C}_k = \ast$, the network $(\calA,\calB,\calC,\calW,\calH)$ reduces to the structured system $(\mathcal{W},\mathcal{H})$. In this special case, Proposition~\ref{p:algebraic} reduces to Proposition~\ref{p:fullrank}. The contribution of Theorem~\ref{t:GGC} and Proposition~\ref{p:algebraic}, however, is that it allows the verification of controllability of a more general class of networked structured systems, where each of the nodes has arbitrary state-space dimension.

\begin{remark}
By Proposition~\ref{p:algebraic}, the analysis of controllability of structured networks boils down to the verification of full row rank of two pattern matrices. Checking whether a given pattern matrix has full rank can be done efficiently, for example, by applying a so-called color change rule to a graph associated to the pattern matrix \cite{Jia2020}.
\end{remark}

We conclude this section with the following corollary of Proposition~\ref{p:algebraic} indicating that  controllability of a structured network requires controllability of the individual node systems.

\begin{corollary}
\label{c:NesCon}
If the network $(\calA, \calB, \calC, \calW, \calH)$ is controllable then $(\mathcal{A}_k,\mathcal{B}_k)$ is controllable for all $k = 1,2,\dots,N$.
\end{corollary}

\begin{IEEEproof}
By Proposition~\ref{p:algebraic}, $(\calA, \calB, \calC, \calW, \calH)$ is controllable if and only if $\begin{bmatrix} \calA + \calB\calW\calC & \calB\calH \end{bmatrix}$ and $\begin{bmatrix} \bar{\calA} + \calB\calW\calC & \calB\calH \end{bmatrix}$ have full row rank. This immediately implies that $\begin{bmatrix} \calA & \calB \end{bmatrix}$ and $\begin{bmatrix} \bar{\calA} & \calB \end{bmatrix}$
have full row rank. By the special structure of $\calA$ and $\calB$ (see \eqref{eq:blkdiagpatterns}), $\begin{bmatrix}\mathcal{A}_k & \calB_k \end{bmatrix}$ and $\begin{bmatrix}\bar{\mathcal{A}}_k & \calB_k \end{bmatrix}$ have full row rank for all $k = 1,2,\dots,N$. By Proposition~\ref{p:fullrank}, this implies that $(\mathcal{A}_k,\mathcal{B}_k)$ is controllable for all $k = 1,2,\dots,N$.
\end{IEEEproof}

\section{Scalable algebraic conditions for controllability of structured networks}\label{s:reduce_algebraic}

In the previous section, we have provided algebraic conditions for controllability of structured networks. These conditions involve checking full row rank of two pattern matrices of dimensions $n \times (n + m)$. However, if the node systems have large state space dimensions $n_k$, then the overall state-space dimension $n = \sum^{N}_{k =1} n_k$ may be prohibitively large.
Therefore, in this section we introduce a new method to verify the full rank property of the pattern matrices in Proposition~\ref{p:algebraic}.
This method will replace these pattern matrices by two auxiliary pattern matrices of much smaller dimensions than $n$.

Specifically, for a given pattern matrix
\begin{equation}
\label{originalpattern}
\bbm \calA + \calB \calW \calC & \calB \calH \ebm,
\end{equation}
we will define a new pattern matrix
\begin{equation}
    \label{newpattern}
\bbm \calhA + \calhB \calW \calhC & \calhB \calH \ebm
\end{equation}
such that \eqref{originalpattern} has full row rank if and only if \eqref{newpattern} has full row rank, but \eqref{newpattern} has much smaller dimensions than \eqref{originalpattern}.

In this new pattern matrix, $\calW$ and $\calH$ remain unchanged, while $(\calA,\calB,\calC)$ is replaced by a reduced system $(\calhA,\calhB,\calhC)$ of the form
\beq \label{eq:pattern2}
\begin{aligned}
\calhA  &= \bdiag(\calhA_{1}, \ldots, \calhA_{N}),\\
\calhB  &= \bdiag(\calhB_{1}, \ldots, \calhB_{N}),\\
\calhC &= \bdiag(\calhC_{1}, \ldots, \calhC_{N}).
\end{aligned}
\eeq
This means that each node system is replaced by a node system $(\calhA_k,\calhB_k,\calhC_k)$ of dimension $\hatn_k$ for $k = 1,2,\ldots,N$.
Obviously, once we have established a procedure to reduce \eqref{originalpattern}, the same procedure can be applied to $\bbm \bar{\calA} + \calB \calW \calC & \calB \calH\ebm$.

In the sequel, we will assume that $\bbm \calA_k & \calB_k \ebm$ has full row rank for $k = 1, \ldots, N$, an assumption that is without loss of generality by Corollary~\ref{c:NesCon}. We will now explain how to define the ``reduced" system $(\calhA, \calhB, \calhC)$ in \eqref{eq:pattern2}. Our strategy will be to replace the node systems $(\mathcal{A}_k,\mathcal{B}_k,\mathcal{C}_k)$ one by one. To this end, we have the following definition.
\bdfn \label{def:replace}
Consider the matrices $\calA, \calB, \calC, \mathcal{W}$ and $\mathcal{H}$, given in \eqref{eqW}, \eqref{eqH} and \eqref{eq:blkdiagpatterns}. Suppose that $(\hat{\calA},\hat{\calB},\hat{\calC})$ is obtained from $(\calA, \calB, \calC)$ by replacing node system $(\calA_k,\calB_k,\calC_k)$ by $(\calhA_k,\calhB_k,\calhC_k)$ for some $k \in \{1, \dots, N\}$. We say that $(\calA_k,\calB_k,\calC_k)$ and $(\calhA_k,\calhB_k,\calhC_k)$ are \emph{equivalent} if \eqref{originalpattern} has full row rank if and only if \eqref{newpattern} has full row rank.
\edfn

The question is now under what conditions two node systems are equivalent.
To answer this question we need the following notion of \emph{independence}.
\bdfn \label{d:independent}
Let $\calM \in \{ 0,\ast,?\}^{1 \times q}$ and $\calN \in \{0,\ast,?\}^{r \times q}$.
We call the pattern vector $\calM$ {\em independent of $\calN$} if for all
$M \in \calP(\calM)$, $N \in \calP(\calN)$ , $z_1 \in \bbR$
and $z_2 \in \bbR^{r}$, $$\bbm z_1 & z_2^\top \ebm \bbm M \\ N \ebm = 0 \text{ implies } z_1 = 0.$$
\edfn

Define $\calA_{k,1}$ as the row in $\calA_{k}$ corresponding to the  position of the $\ast$ entry in $\calB_k$ and $\calA_{k,2}$ as the pattern matrix obtained from $\calA_k$ by removing the row $\calA_{k,1}$.
By our standing hypothesis that $\bbm \calA_k & \calB_k \ebm$ has full row rank, it is clear that $\calA_{k,2}$ has full row rank.
 We now distinguish the following four properties of $(\calA_{k},\calB_k,\calC_k)$:
\begin{enumerate}[(S1)]
\item \label{S1} The pattern vector $\calA_{k,1}$ is independent of $\col(\calA_{k,2},\calC_{k})$.
\item The pattern vector $\calA_{k,1}$ is independent of $ \calA_{k,2}$. \label{S2}
\item The pattern vector $\calC_{k}$ is independent of $\calA_{k}$. \label{S3}
\item The pattern vector $\calC_{k}$ is independent of $\calA_{k,2}$. \label{S4}
\end{enumerate}
 Clearly, \ref{S1} implies \ref{S2}, and \ref{S3} implies \ref{S4}.
Note that \ref{S2} is equivalent to saying that $\mathcal{A}_k$ has full row rank.
 Hence, \ref{S2} and \ref{S3} are mutually exclusive.
 Moreover, it also holds that \ref{S1} and \ref{S4} are mutually exclusive.
 Therefore, for any given $(\calA_k, \calB_k, \calC_k)$ {\em exactly one} of the following six conditions holds:
\begin{enumerate}[(C1)]
\item Property \ref{S1} holds. \label{C1}
\item Property \ref{S3} holds. \label{C2}
\item Properties \ref{S2} and \ref{S4} hold. \label{C3}
\item Property \ref{S2} holds but neither \ref{S1} nor \ref{S4} holds. \label{C4}
\item Property \ref{S4} holds but neither \ref{S2} nor \ref{S3} holds. \label{C5}
\item Neither \ref{S2} nor \ref{S4} holds. \label{C6}
\end{enumerate}

The following main result of this section now gives necessary and sufficient conditions under which two node systems are equivalent in the sense of Definition \ref{def:replace}.

\bthm \label{t:AC4equal}
Consider two node systems $(\calA_k, \calB_k, \calC_k)$ and $(\calhA_k,\calhB_k,\calhC_k)$ both having state-space dimension at least two. Then $(\calA_k,\calB_k,\calC_k)$ and $(\calhA_k,\calhB_k,\calhC_k)$ are equivalent if and only if they satisfy the same condition \Ci{}.
\ethm

The power of Theorem~\ref{t:AC4equal} becomes clear once we realize that for each $i = 1,2,\ldots,6$ there exists a node system of dimension  $2$ that satisfies condition \Ci{}. As a consequence, any system $(\calA_k,\calB_k,\calC_k)$ of \emph{arbitrary dimension} can be reduced to a system $(\calhA_k,\calhB_k,\calhC_k)$ of dimension at most $2$. It turns out that for systems that satisfy \ref{C2}, \ref{C3} and \ref{C5} we can even reduce $(\calA_k,\calB_k,\calC_k)$ to a scalar node system.
In particular, for each $i = 1,2, \ldots, 6$ we define a {\em standard node system} $(\calhA_{\Ci},\calhB_{\Ci},\calhC_{\Ci})$, where
\beq \label{eq:standard_A}
\begin{aligned}
&\calhA_{\textup{(C1)}} = \bbm 0 & \ast\\ \ast & 0 \ebm, \quad \calhA_{\textup{(C2)}} = 0, \quad \calhA_{\textup{(C3)}} = \ast,\\
&\calhA_{\textup{(C4)}} = \bbm 0 & \ast\\ \ast & ? \ebm, \quad \calhA_{\textup{(C5)}} = ?, \quad \calhA_{\textup{(C6)}} = \bbm 0 & 0 \\ \ast & 0 \ebm,\\
\end{aligned}
\eeq
and
\beq \label{eq:standard_BC}
\calhB_{\Ci}^\top = \calhC_{\Ci} : =
\begin{cases}
\vspace{10pt}
 \bbm \ast & 0\ebm & \quad i = 1, 4, 6\\
~~~~\ast & \quad i = 2, 3, 5.
\end{cases}
\eeq

\blem \label{l:candidate_1}
Suppose that $(\calA_k, \calB_k, \calC_k)$ satisfies \Ci{}.
Then $(\calA_k, \calB_k, \calC_k)$ and $(\calhA_{\Ci},\calhB_{\Ci},\calhC_{\Ci})$ (as defined by \eqref{eq:standard_A}, \eqref{eq:standard_BC}) are equivalent.
\elem

In that case, the node system $(\calhA_{\Ci},\calhB_{\Ci},\calhC_{\Ci})$ is called {\em the standard node system corresponding to} $(\calA_k, \calB_k, \calC_k)$.

To sum up the results from this section, we state the following theorem that provides a new algebraic condition for controllability of the structured network $(\calA,\calB,\calC,\calW,\calH)$.
\bthm \label{t:ACSC}
Consider the network $(\calA,\calB,\calC,\calW,\calH)$.
Define $(\calhA,\calhB,\calhC)$ as in \eqref{eq:pattern2}, where $(\calhA_k,\calhB_k,\calhC_k)$ is the standard node system corresponding to $(\calA_k, \calB_k, \calC_k)$ for $k = 1, \ldots, N$.
Similarly, define $(\hat{\bar{\calA}}, \hat{\bar{\calB}},\hat{\bar{\calC}})$ as
\begin{equation}
\label{eq:blkdiagpatterns_scalable_2}
\begin{aligned}
 	 \calhbA &= \bdiag(\calhbA_1, \ldots, \calhbA_N),\\
	\calhbB &= \bdiag(\calhbB_1, \ldots, \calhbB_N),\\
	\calhbC &= \bdiag(\calhbC_1, \ldots, \calhbC_N)
\end{aligned}
\end{equation}
where $(\calhbA_k, \calhbB_k, \calhbC_k)$ is the standard node system corresponding to $(\calbA_k, \calB_k, \calC_k)$ for $k = 1, \ldots, N$.
Then $(\calA,\calB,\calC,\calW,\calH)$ is controllable if and only if both
$\bbm \calhA + \calhB \calW\calhC & \calhB \calH \ebm$
 and $\bbm \hat{\bar{\calA}} + \hat{\bar{\calB}} \calW\hat{\bar{\calC}} & \hat{\bar{\calB}} \calH \ebm$ have full row rank.
\ethm

\bre \label{r:homogeneous}
Note that the pattern matrices $\bbm \calhA + \calhB \calW\calhC & \calhB \calH\ebm$ and $\bbm \hat{\bar{\calA}} + \hat{\bar{\calB}} \calW\hat{\bar{\calC}} & \hat{\bar{\calB}} \calH \ebm$ appearing in Theorem~\ref{t:ACSC} have at most $2N$ rows.
In general, these pattern matrices are thus of much lower dimension than the dimension $n = \sum_{k = 1}^N n_k$ of the original network. Controllability analysis becomes particularly simple in the case that $(\calA,\calB,\calC,\calW,\calH)$ is \emph{homogeneous} (referred to as {\em similar} in  \cite{CPAKJ2017}), i.e., if all node systems $(\calA_k,\calB_k,\calC_k)$ are identical.
Indeed, in this scenario we only need to check which of the $6$ conditions \ref{C1} to \ref{C6} are satisfied for $(\calA_1, \calB_1, \calC_1)$ and $(\bar{\calA}_1, \calB_1, \calC_1)$.
Subsequently, the systems $(\calA,\calB,\calC)$ and $(\calbA,\calB,\calC)$ can be reduced to \eqref{eq:pattern2} and \eqref{eq:blkdiagpatterns_scalable_2} respectively, where the reduced patterns are also ``homogeneous" in the sense that, e.g., $\calhA_1 = \cdots = \calhA_N$, etc.
\ere

\bex \label{ex:2}
We revisit the structured network $(\calA,\calB,\calC,\calW,\calH)$ of Example~\ref{ex:1}. Our aim will be to reduce each of the node systems $(\calA_k,\calB_k,\calC_k)$ and $(\bar{\calA_k},\calB_k,\calC_k)$ ($k = 1,\dots,7$), so that we can verify controllability of $(\calA,\calB,\calC,\calW,\calH)$ by assessing the rank of lower dimensional pattern matrices. We start with $(\mathcal{A}_1,\calB_1,\calC_1)$ which was given by
\begin{equation}
    \label{A1B1C1}
\begin{aligned}
\calA_1 &=
\bbm
0 & \ast & 0 & 0\\
\ast & \ast & \ast & 0\\
0& 0 &0 & \ast\\
\ast & 0 &0 & 0\\
\ebm, \quad \calB_1 = \bbm 0 \\ \ast \\ 0 \\ 0\ebm, \quad \calC_1^\top = \bbm \ast \\ 0 \\ 0 \\ 0 \ebm.
\end{aligned}
\end{equation}
Note that $\calA_{1,1} = \begin{bmatrix} \ast & \ast & \ast & 0 \end{bmatrix}$ is independent of the matrix
$$
\col(\calA_{1,2}, \calC_1) = \begin{bmatrix}
0 & \ast & 0 & 0\\
0& 0 &0 & \ast\\
\ast & 0 &0 & 0\\
\ast & 0 & 0 & 0
\end{bmatrix}.
$$
That is, property \ref{S1} holds.
Hence, $(\mathcal{A}_1,\calB_1,\calC_1)$ satisfies condition \ref{C1}, and consequently, its standard node system is given by
\begin{align*}
\calhA_1 &=
\bbm
0 & \ast \\
\ast & 0
\ebm, \quad \calhB_1 = \calhC_1^\top = \bbm \ast \\ 0 \ebm.
\end{align*}
In a similar manner, we can reduce the other node systems. This results in the lower dimensional pattern matrices $\bbm \calhA + \calhB \calW \calhC & \calhB \calH \ebm$ and $\bbm \calhbA + \calhbB \calW \calhbC & \calhbB \calH \ebm$. We provide a graphical visualization of the original pattern matrices, as well as their reduced counterparts in Figures~\ref{f:ABWC} and \ref{f:re_Net}. It can be verified that the reduced pattern matrices have full row rank. This can, for instance, be done by applying the color change rule \cite{Jia2020} to the reduced graphs in Figures~\ref{f:rABCWH} and \ref{f:rbABCWH}. Therefore, by Theorem~\ref{t:ACSC} we conclude that the structured network $(\calA,\calB,\calC,\calW,\calH)$ is controllable.

\input{f_ABWC}

\input{f_re_Net}

\eex

\section{Graph theoretic conditions}\label{s:reduce_graph}

To apply Theorem~\ref{t:ACSC}, we need to know which condition \Ci{} is satisfied for each node system $(\calA_{k}, \calB_{k}, \calC_{k})$. This means that we have to check which of the properties \ref{S1}-\ref{S4} are satisfied for each of the node systems. Note that these four properties all involve the independence of certain pattern vectors. Although full row rank of pattern matrices can be checked efficiently \cite[Theorem 10]{Jia2020}, we are not aware of any methods to check whether a pattern vector is independent of a pattern matrix. Therefore, in this section we provide a graph theoretic method to verify whether a given pattern vector is independent of a pattern matrix.

Before we explain the procedure, we recall some graph theoretic preliminaries from \cite{Jia2020}. Define the directed graph associated with $\calM \in \{0,\ast,?\}^{p \times q}$ as $G(\calM) = (V,E)$. Here the node set $V$ is given by $V = \{1, \ldots, \max(p,q)\}$ and the edge set $E$ is defined as
$$
E = \{(i,j) \in V\times V \mid \calM_{ji} = \ast \mbox{ or } ?\}.
$$
To distinguish between $\ast$ and $?$ entries in $\calM$, we partition the edge set $E$ into two disjoint subsets $E_\ast$ and $E_?$ given by
$$E_\ast = \{(i,j) \in E \mid \calM_{ji} = \ast\}, \: E_? = \{(i,j) \in E \mid \calM_{ji} = ?\}.
$$
Consider the following coloring procedure which was defined in \cite{Jia2020}:
\begin{enumerate}
	\item Initially, color all nodes of $G(\calM)$ white.
	\item If a node $i$ has exactly one white out-neigbor $j$ and $(i,j) \in E_\ast$, change the color of $j$ to black. \label{i:step2}
	\item Repeat step~\ref{i:step2} until no more changes are possible.
\end{enumerate}
The \emph{derived set} $S(\calM)$ of $G(\calM)$ is defined as the set of all black nodes obtained by applying the above procedure to $G(\calM)$.
It has been shown in \cite[Theorem 10]{Jia2020} that $\mathcal{M}$ has full row rank if and only if $S(\calM) =  \nset{p}$. In what follows, we use a similar idea to give graph theoretic conditions under which a pattern vector is independent of a pattern matrix.

\blem \label{l:in_color}
Let $\calM \in \{0,\ast,?\}^{1 \times p}$ and $\calN \in \{0,\ast,?\}^{r \times p}$. Consider the graph $G(\col(\calM,\calN)) = (V,E)$. Then $\calM$ is independent of $\calN$ if and only if node $1 \in V$ is contained in the derived set $S(\col(\calM,\calN))$.
\elem

Lemma~\ref{l:in_color} can be immediately applied to check which of the mutually exclusive conditions \ref{C1}-\ref{C6} holds for a given node system. We illustrate this in the following example.

\bex \label{ex:3}
We revisit the network in Examples \ref{ex:1} and \ref{ex:2}. The purpose of this example is to apply the graph theoretic test of Lemma~\ref{l:in_color} to show that the system $(\calA_1,\calB_1,\calC_1)$ in \eqref{A1B1C1} satisfies condition \ref{C1}.
To do so, consider the graph $G(\col(\calA_1, \calC_1))$ depicted in Figure~\ref{f:A_1_1}. Initially, color all nodes in this graph white. Clearly, node $3$ has only one white out-neighbor $2$, and $(3,2) \in E_\ast$. We thus color $2$ black. Similarly, $3$ is colored by $4$.
Finally, node $1$ is colored by $2$. No more nodes can be colored,  hence the coloring process stops and we obtain the derived set $S(\col(\calA_1, \calC_1)) = \{1,2,3\}$. Since the second entry of $\calB_1$ is equal to $\ast$ and $2 \in S(\col(\calA_1, \calC_1))$, we see that property \ref{S1} holds. Therefore, condition \ref{C1} is satisfied for $(\calA_1,\calB_1,\calC_1)$.

\begin{figure}[h!]
\centering
    \begin{tikzpicture}[->,>=stealth',shorten >=1pt,auto,node distance=2.8cm, semithick]
		\tikzset{VS1/.style = {shape = circle,
                color=black,
                fill=white!80!black,
                minimum size=0.5cm,
                text = black,
                inner sep = 2pt,
                outer sep = 1pt,
                draw}
        }

        \foreach \theta in {1,2,3,4}{\node[VS1] (1\theta) at (\theta*90: 1cm) {$\theta$};}
        \node[VS1] (15) at (2,1) {$5$};

        \draw (11) to [out=-145, in=55, looseness = 1] (12);
        \draw (11) to (14);
        \draw (12) to [out=35, in=-125, looseness = 1] (11);
        \draw (12) to [out=-115, in=-155, looseness = 10] (12);
        \draw (14) to (13);
        \draw (13) to (12);

        \draw (11) to (15);
	\end{tikzpicture}
	\caption{The graph $G(\col(\calA_1,\calC_1))$.}
    \vspace{-3mm}
	\label{f:A_1_1}
\end{figure}

\eex

\bre \label{r:gt_cases}
Suppose that we want to decide which one of the conditions \ref{C1}-\ref{C6} holds for a given node system. In the worst case, we have to check each of the four properties \ref{S1}-\ref{S4} one by one. By Lemma~\ref{l:in_color}, this boils down to computing the derived set of three different (but strongly related) graphs. It turns out that it is not necessary to recompute the entire derived set in each of these graphs. In fact, in this remark we provide a more efficient procedure to check which one of the conditions \ref{C1}-\ref{C6} holds.

Let $j_k$ be such that $j_k$-th entry of $\calB_{k}$ is equal to $\ast$.
Let $\calT = \col(\calA_k,\calC_k)$ and apply the color change rule to $G(\calT)$ in order to compute the derived set $S(\calT)$.
\begin{itemize}
    \item $(\calA_k,\calB_k,\calC_k)$ satisfies \ref{S1} if and only if the vertex $j_k$ is contained in $S(\calT)$. If $j_k \in S(\calT)$ then we are done since we know that $(\calA_k,\calB_k,\calC_k)$ satisfies condition \ref{C1}.
\item $(\calA_k,\calB_k,\calC_k)$ satisfies \ref{S3} if and only if the vertex $n_k + 1$ is contained in $S(\calT)$. Again, if $n_k+1 \in S(\calT)$ we are done, as in this case $(\calA_k,\calB_k,\calC_k)$ satisfies condition \ref{C2}.

If neither \ref{S1} nor \ref{S3} holds, we move on to check whether properties \ref{S2} and/or \ref{S4} hold.
\item Color the vertices in $S(\calT) \cup \{n_k + 1\}$ black.
Apply the color change rule on $G(\calT)$ until no more color changes are possible, and let $S'$ be the resulting set of black vertices.
Then property \ref{S2} holds if and only if $j_k \in S'$.
\item Color the vertices in $S(\calT) \cup \{j_k\}$ black.
Apply the color change rule on $G(\calT)$ until no more color changes are possible, and let $S''$ be the resulting set of black nodes.
It can then be shown that property \ref{S4} holds if and only if $n_k + 1 \in S''$.
\item Depending on which of the properties \ref{S2} and \ref{S4} hold, we can easily determine which condition \ref{C3}-\ref{C6} is satisfied for $(\calA_k,\calB_k,\calC_k)$.
\end{itemize}
The above remark outlines a conceptual algorithm to check which condition \ref{C1}-\ref{C6} is satisfied for a given node system. We again emphasize that this procedure is efficient in the sense that it avoids recomputing the entire derived sets of the different graphs related to \ref{S1}-\ref{S4}.
\ere

\section{Conclusion and Discussion}\label{s:s5_con}

In this paper, we have studied strong structural controllability of structured networks.
In contrast to existing work, where the node systems are usually assumed to be single integrators, in this paper, we allow single-input single-output node systems with arbitrary state space dimensions.
The node systems and the structured interconnection laws interconnecting these are all represented by pattern matrices with three possible entries, namely $0$, nonzero indeterminate, and arbitrary (zero or nonzero) indeterminate.
We have proven that a structured network is controllable if and only if an associated structured system is controllable.
This makes it possible to check controllability of a structured network by applying existing tests \cite{Jia2020} for controllability of structured systems. Applying these existing tests might be intractable because of the large state space dimension of the structured network. In order to overcome this difficulty, we have
shown that controllability of a given structured network can be tested by replacing the original network by a new network in which all original node systems have been replaced by (auxiliary) node systems from a set of six standard node systems with state space dimensions either 1 or 2. This means that controllability of any network can be checked by testing controllability of a structured system of state space dimension at most twice the number of node systems, regardless of the dimensions of the original node systems. As such, this method is scalable, because after replacing in the network one of the original node systems by a possibly higher dimensional node system, testing controllability will only involve a check which of the six standard (first or second order) standard node systems should be used as its substitute. In order to determine which of the six standard node systems should replace a given original node system, we have introduced a color change procedure to be applied to the graph of each original node system.

We conclude this section with some suggestions for future research. Whereas the present paper deals only with single-input single-output node systems, a venue for future research could be to generalize our results to general multi-input multi-output node systems.
Another opportunity for future research is to extend  our results to a wider range of system properties.
Obviously, tests for structural observability of networks can be obtained by dualizing our results.
Other structural properties of interest are, for example, input-state observability and output controllability \cite{Shali2020}, fault detection and isolation \magenta{\cite{JTC2020-FDI}}, and system invertibility \magenta{\cite{TSH2012}} of structured networks.

\appendix[]
\section{Proofs}\label{s:Proofs}
In this Appendix, for given pattern matrices $\calM_1$ and $\calM_2$,
we will denote the Cartesian product $\calP(\calM_1) \times \calP(\calM_2)$ by $\calP(\calM_1, \calM_2)$, and likewise for three or more pattern matrices.
\subsection{Proof of Theorem \ref{t:AC4equal}}
For the proof of Theorem \ref{t:AC4equal}, the following auxiliary result will be instrumental:
\blem\label{l:condition}
Let $k \in \{1,\dots,N\}$. Suppose that for all $(A_k, B_k, C_k) \in \calP(\calA_k, \calB_k, \calC_k)$, $x_k\in\bbR^{n_k}$ and $\lambda_k\in\bbR$ such that
	\begin{equation}\label{eq:conditions_from}
	x_k\neq 0,\quad x_k^\top A_k = \lambda_k C_k,
	\end{equation}
	there exist $(\hatA_k, \hatB_k, \hatC_k)\in\calP(\calhA_k,\calhB_k, \calhC_k)$ and $y_k\in\bbR^{\hatn_k}$ such that
	\begin{equation}\label{eq:conditions_to}
	y_k \neq 0,\quad y_k^\top \hatA_k = \lambda_k \hatC_k, \quad y_k^\top\hatB_k = x_k^\top B_k.
	\end{equation}
	Then $\bbm \calA+\calB\calW\calC & \calB\calH \ebm$ has full row rank if $\bbm \calhA+\calhB\calW\calhC & \calhB\calH \ebm$ has full row rank, where $(\calhA_i,\calhB_i,\calhC_i) = (\calA_i,\calB_i,\calC_i)$ for $i \neq k$.
\elem
\begin{IEEEproof} 
	Suppose that $\bbm \calhA + \calhB \calW\calhC & \calhB \calH \ebm$ has full row rank, but that, on the contrary, $\bbm \calA + \calB\calW\calC & \calB\calH \ebm$ does not have full row rank. Then there exists $(A,B,C,W,H)\in\calP(\calA, \calB, \calC, \calW, \calH) \text{ and nonzero } x \in\bbR^{n}$ such that
	\begin{equation}\label{eq:big_x}
		x^\top \bbm A + BWC & BH \ebm = 0.
	\end{equation}
	Partition $x= \col(x_1, \ldots, x_N)$, where $x_i \in \bbR^{n_i}$ for $i=1,\dots,N$. Then \eqref{eq:big_x} can be rewritten as
	\begin{equation}\label{eq:small_x1}
	\sum_{j=1}^N x_j ^\top B_j H_j = 0 
	\end{equation}
and
	\begin{equation}\label{eq:small_x2}
	 x_i^\top A_i + \lr*{\sum_{j=1}^N w_{ji} x_j^\top B_j} C_i = 0 
	\end{equation}
	for $i=1,\dots,N$.
	If $x_k = 0$, then take $y_k = 0$ and arbitrary $(\hatA_k, \hatB_k, \hatC_k)\in\calP(\calhA_k, \calhB_k, \calhC_k)$. Otherwise, $x_k \neq 0$ and $\lambda_k = -\sum_{j=1}^N w_{jk}x_j^\top B_j$ are such that \eqref{eq:conditions_from} holds, hence there exist $(\hatA_k, \hatB_k, \hatC_k)\in\calP(\calhA_k, \calhB_k, \calhC_k)$ and $y_k\in\bbR^{\hatn_k}$ such that \eqref{eq:conditions_to} holds. 
Therefore, in both cases we have that $$y_k^\top \hatA_k = \lambda_k \hatC_k \qand y_k^\top \hat B_k = x_k^\top B_k,$$ while $y_k = 0$ if and only if $x_k = 0$. 
Let $y_i = x_i$, $\hat A_i = A_i$, $\hat B_i = B_i$ and $\hat C_i = C_i$ for all $i\neq k$. Then $y_j^\top \hat B_j = x_j^\top B_j$ for all $j=1,\dots,N$, and thus \eqref{eq:small_x1} and \eqref{eq:small_x2} imply that
		\begin{equation}\label{eq:small_y_no_k1}
		\sum_{j=1}^N y_j^\top \hat B_j H_j = 0
		\end{equation}
and
			\begin{equation}\label{eq:small_y_no_k2}
y_i^\top \hat A_i + \lr*{\sum_{j=1}^N w_{ji}y_j^\top \hat B_j} \hat C_i = 0
		\end{equation}
		for all $i \neq k$.
		Furthermore, $\lambda_k = -\sum_{j=1}^N w_{jk}y_j^\top \hat B_j$ and $y_k^\top \hat A_k = \lambda_k \hat C_k$ imply that \eqref{eq:small_y_no_k2} holds for all $i = 1,\dots, N$.

		 Let $y = \col(y_1, \ldots, y_N)$,  $\hatA = \bdiag(\hat A_1, \ldots, \hat A_N)$, $\hat B = \bdiag(\hat B_1, \ldots, \hat B_N)$ and $\hat C = \bdiag(\hat C_1, \ldots, \hat C_N)$. 
		It then follows that
		  $(\hatA, \hatB, \hatC)\in\calP(\calhA, \calhB, \calhC)$, and
	\begin{equation}\label{eq:big_y}
	 y \neq 0, \quad y^\top \bbm \hat A + \hat BW \hat C & \hat BH \ebm = 0.
	\end{equation}
	Therefore, we reach a contradiction, and hence $\bbm \calA + \calB\calW\calC & \calB\calH \ebm$ has full row rank.
\end{IEEEproof}

We are now ready to provide a Proof of Theorem \ref{t:AC4equal}.

\begin{IEEEproof}[Proof of Theorem \ref{t:AC4equal}]
	Suppose that $(\calA_k,\calB_k,\calC_k)$ and $(\calhA_k,\calhB_k,\calhC_k)$ satisfy the same condition \Ci{}. 
	Recall that $(\calA_k,\calB_k,\calC_k)$ and $(\calhA_k,\calhB_k,\calhC_k)$ are equivalent if $\bbm \calA+\calB\calW\calC & \calB\calH \ebm$ has full row rank if and only if $\bbm \calhA+\calhB\calW\calhC & \calhB\calH \ebm$ has full row rank, where $(\calhA_i,\calhB_i,\calhC_i) = (\calA_i,\calB_i,\calC_i)$ for $i \neq k$. 
	We will first show that  $\bbm \calA+\calB\calW\calC & \calB\calH \ebm$ has full row rank if $\bbm \calhA + \calhB\calW\calhC & \calhB\calH \ebm$ has full row rank. Recall from Section \ref{s:reduce_algebraic} that $\calA_{k,1}$ is defined as the row in $\calA_{k}$ corresponding to the  position of the $\ast$ entry in $\calB_k$ and $\calA_{k,2}$ as the pattern matrix obtained from $\calA_k$ by removing the row $\calA_{k,1}$.
	In view of Lemma~\ref{l:condition}, it suffices to show that for all $(A_{k,1}, A_{k,2}, C_k)\in\calP(\calA_{k,1}, \calA_{k,2}, \calC_k)$, $x_{k,1} \in \bbR$, $x_{k,2} \in\bbR^{n_k - 1}$ and $\lambda_k\in\bbR$ such that 
	\begin{equation}\label{eq:condition_if}
	\bbm x_{k,1} \\ x_{k,2} \ebm \neq 0 \qand \bbm x_{k,1} & x_{k,2}^\top \ebm \bbm A_{k,1} \\ A_{k,2} \ebm = \lambda_k C_k,
	\end{equation}
	there exist $(\hat A_{k,1}, \hat A_{k,2}, \hat C_k) \in\calP(\calhA_{k,1}, \calhA_{k,2}, \calhC_k)$ and $y_{k,2} \in\bbR^{\hat n_k - 1}$ such that 
	\begin{equation}\label{eq:condition_then}
	\bbm x_{k,1} \\ y_{k,2} \ebm \neq 0 \qand \bbm x_{k,1} & y_{k,2} ^\top \ebm \bbm \hat A_{k,1} \\ \hat A_{k,2}\ebm = \lambda_{k} \hat C_{k}.
	\end{equation}
	 With this in mind, let $(A_{k,1}, A_{k,2}, C_k)\in\calP(\calA_{k,1}, \calA_{k,2}, \calC_k)$, $x_{k,1} \in \bbR$, $x_{k,2} \in\bbR^{n_k - 1}$ and $\lambda_k\in\bbR$ be such that \eqref{eq:condition_if} holds. Note that $x_{k,1}$ and $\lambda_{k}$ are not both zero. 
	 Indeed, if $x_{k,1} = \lambda_k = 0$, then $x_{k,2} = 0$ since $\calA_{k,2}$ is assumed to have full row rank and we reach a contradiction. 
	 We will consider each of the conditions \ref{C1}, \ldots, \ref{C6} separately.
	
	\textbf{Condition~\ref{C1}.} Since \ref{S1}  holds for $(\calA_{k,1}, \calA_{k,2}, \calC_k)$, \eqref{eq:condition_if} implies that $x_{k,1} = 0$, hence $\lambda_k \neq 0$.
	Furthermore, since \ref{S1}  holds for $(\calhA_{k,1}, \calhA_{k,2}, \calhC_k)$, it follows that \ref{S4}  does not hold, i.e., $\calhC_k$ is not independent of $\calhA_{k,2}$. 
	Given that $\lambda_k \neq 0$, this implies that there exist $( \hat A_{k,2}, \hat C_{k}) \in\calP(\calhA_{k,2}, \calhC_k)$ and nonzero $y_{k,2} \in\bbR^{\hat n_k - 1}$ such that $y^\top_{k,2} \hat A_{k,2} = \lambda_k\hat C_{k}$. Then $\eqref{eq:condition_then}$ is satisfied for all $\hat A_{k,1}\in\calP(\calhA_{k,1})$.
	
	\textbf{Condition~\ref{C2}.} Since \ref{S3}  holds for $(\calA_{k,1}, \calA_{k,2}, \calC_k)$, \eqref{eq:condition_if} implies that $\lambda_k = 0$, hence $x_{k,1}\neq 0$. Furthermore, since \ref{S3}  holds for $(\calhA_{k,1}, \calhA_{k,2}, \calhC_k)$, it follows that \ref{S2}  does not hold and $\calhA_{k,1}$ is not independent of $\calhA_{k,2}$. Given that $x_{k,1} \neq 0$, this implies that there exist $( \hatA_{k,1}, \hatA_{k,2}) \in\calP(\calhA_{k,1}, \calhA_{k,2})$ and nonzero $y_{k,2} \in\bbR^{\hat n_k - 1}$ such that $y_{k,2}^\top \hatA_{k,2} = -x_{k,1}\hat A_{k,1}$. 
	Then $\eqref{eq:condition_then}$ is satisfied for all $\hat C_k\in\calP(\calhC_k)$.
	
	\textbf{Condition~\ref{C3}.} We claim that since \ref{S2}  and \ref{S4}  hold for $(\calA_{k,1}, \calA_{k,2}, \calC_k)$, \eqref{eq:condition_if} implies that $\lambda_k \neq 0$ and $x_{k,1} \neq 0$. Indeed, if $\lambda_k = 0$, then \ref{S2}  implies that $x_{k,1} = 0$, and if $x_{k,1} = 0$, then \ref{S4}  implies that $\lambda_k = 0$. But we know that $\lambda_k$ and $x_{k,1}$ are not both zero, hence we reach a contradiction. Since \ref{S1}  does not hold for $(\calhA_{k,1}, \calhA_{k,2}, \calhC_k)$, this implies that there exist $(\hat A_{k,1}, \hat A_{k,2}, \hat C_k')\in\calP(\calhA_{k,1}, \calhA_{k,2}, \calhC_k)$, $y_{k,2}\in\bbR^{\hat n_k - 1}$ and $\mu_k\in\bbR$ such that
	\begin{equation*}
	\bbm x_{k,1} & y_{k,2}^\top \ebm \bbm \hat A_{k,1} \\ A_{k,2} \ebm = \mu_k \hat C_k',
	\end{equation*}
	But then $\mu_k \neq 0$ because \ref{S2}  holds for $(\calhA_{k,1}, \calhA_{k,2}, \calhC_k)$, hence $\hat C_k = \tfrac{\mu_k}{\lambda_k}\hat C_k'\in\calP(\calhC_k)$ is such that \eqref{eq:condition_then} holds.
	
	\textbf{Condition~\ref{C4}.} In condition \ref{C3} we showed that $\lambda_k \neq 0$ whenever \ref{S2}  holds for $(\calA_{k,1}, \calA_{k,2}, \calC_k)$. With this in mind, suppose that $x_{k,1} = 0$. Given that \ref{S4}  does not hold for $(\calhA_{k,1}, \calhA_{k,2}, \calhC_k)$ and $\lambda_k \neq 0$, there exist $(\hat A_{k,2}, \hat C_k)\in\calP(\calhA_{k,2}, \calhC_k)$ and nonzero $y_{k,2}\in\bbR^{\hat n_k - 1}$ such that $y_{k,2}^\top \hat A_{k,2} = \lambda_k \hat C_k$, hence \eqref{eq:condition_then} is satisfied for all $\hat A_{k,1}\in\calP(\calhA_{k,1})$. Conversely, suppose that $x_{k,1} \neq 0$. Since \ref{S1}  does not hold for $(\calhA_{k,1}, \calhA_{k,2}, \calhC_k)$, this implies that there exist $(\hat A_{k,1}, \hat A_{k,2}, \hat C_k')\in\calP(\calhA_{k,1}, \calhA_{k,2}, \calhC_k)$, $y_{k,2}\in\bbR^{\hat n_k -1}$ and $\mu_k\in\bbR$ such that
	\begin{equation*}
	\bbm x_{k,1} & y_{k,2}^\top \ebm \bbm \hat A_{k,1} \\ \hat A_{k,2} \ebm = \mu_k \hat C_k'.
	\end{equation*}
	But then $\mu_k\neq 0$ because \ref{S2}  holds for $(\calhA_{k,1}, \calhA_{k,2}, \calhC_k)$.
	 Hence, $\hat C_k = \tfrac{\mu_k}{\lambda_k}\hat C_k'\in\calP(\calhC_k)$ is such that \eqref{eq:condition_then} holds.
	
	\textbf{Condition~\ref{C5}.}  In condition \ref{C3} we showed that $x_{k,1}\neq 0$ whenever \ref{S4}  holds for $(\calA_{k,1}, \calA_{k,2}, \calC_k)$. With this in mind, suppose that $\lambda_k = 0$. Given that \ref{S2}  does not hold for $(\calhA_{k,1}, \calhA_{k,2}, \calhC_k)$ and $x_{k,1} \neq 0$, there exist $(\hat A_{k,1}, \hat A_{k,2})\in \calP(\calhA_{k,1}, \calhA_{k,2})$ and $y_{k,2}\in\bbR^{\hat n_k -1}$ such that $y_{k,2}^\top \hat A_{k,2} = -x_{k,1}\hat A_{k,1}$, hence \eqref{eq:condition_then} is satisfied for all $\hat C_k\in\calP(\calhC_k)$. Conversely, suppose that $\lambda_k\neq 0$. Since \ref{S3}  does not hold for $(\calhA_{k,1}, \calhA_{k,2}, \calhC_k)$, this implies that there exist $(\hat A_{k,1}', \hat A_{k,2}, \hat C_k)\in\calP(\calhA_{k,1}, \calhA_{k,2}, \calhC_k)$, $y_{k,1}\in\bbR$ and $y_{k,2}\in\bbR^{\hat n_k -1}$ such that
	\begin{equation*}
	\bbm y_{k,1} & y_{k,2}^\top \ebm \bbm \hat A_{k,1}' \\ \hat A_{k,2} \ebm = \lambda_k \hat C_k.
	\end{equation*}
	But then $y_{k,1} \neq 0$ because \ref{S4}  holds for $(\calhA_{k,1}, \calhA_{k,2}, \calhC_k)$, hence $\hat A_{k,1} = \tfrac{y_{k,1}}{x_{k,1}}\hat A_{k,1}'\in\calP(\calhA_{k,1})$ is such that \eqref{eq:condition_then} holds.
	
	\textbf{Condition~\ref{C6}.} Since $\lambda_k$ and $x_{k,1}$ are not both zero, there are only three cases to consider: $\lambda_k =0$ and $x_{k,1}\neq 0$; $\lambda_k \neq 0$ and $x_{k,1} = 0$; $\lambda_k\neq 0$ and $x_{k,1}\neq 0$. 
	
	To begin with, suppose that $\lambda_k =0$ and $x_{k,1}\neq 0$. Given that \ref{S2}  does not hold for $(\calhA_{k,1}, \calhA_{k,2}, \calhC_k)$, there exist $(\hat A_{k,1}, \hat A_{k,2})\in\calP(\calhA_{k,1}, \calhA_{k,2})$ and $y_{k,2}\in\bbR^{\hat n_k -1}$ such that $y_{k,2}^\top \hat A_{k,2} = -x_{k,1}\hat A_{k,1}$, hence \eqref{eq:condition_then} is satisfied for all $\hat C_k\in\calP(\calhC_k)$. 
	
	Next, suppose that $\lambda_k \neq 0$ and $x_{k,1} = 0$. Given that \ref{S4}  does not hold for $(\calhA_{k,1}, \calhA_{k,2}, \calhC_k)$ and $\lambda_k \neq 0$, there exist $(\hat A_{k,2}, \hat C_k)\in\calP(\calhA_{k,2}, \calhC_k)$ and nonzero $y_{k,2}\in\bbR^{\hat n_k - 1}$ such that $y_{k,2}^\top \hat A_{k,2} = \lambda_k \hat C_k$, hence \eqref{eq:condition_then} is satisfied for all $\hat A_{k,1}\in\calP(\calhA_{k,1})$. 
	
	Finally, suppose that $\lambda_k\neq 0$ and $x_{k,1}\neq 0$. 
	We will distinguish two cases depending on whether there exist $(\hat A_{k,1}, \hat A_{k,2}, \hat C_k)\in\calP(\calhA_{k,1}, \calhA_{k,2}, \calhC_k)$ such that neither $\hat A_{k,1}$ nor $\hat C_k$ is independent of $\hat A_{k,2}$. 
	First, suppose that such $(\hat A_{k,1}, \hat A_{k,2}, \hat C_k)\in\calP(\calhA_{k,1}, \calhA_{k,2}, \calhC_k)$ exists. 
	Then there exist $y_{k,2}',y_{k,2}''\in\bbR^{\hat n_k - 1}$ such that \[y_{k,2}^{\prime\top}A_{k,2}=-x_{k,1}A_{k,1} \qand y_{k,2}^{\prime\prime\top}A_{k,2}=\lambda_kC_{k},\] hence $y_{k,2} = y_{k,2}' + y_{k,2}''$ is such that \eqref{eq:condition_then} holds. 
	Conversely, suppose that all $(\hat A_{k,1}, \hat A_{k,2}, \hat C_k)\in\calP(\calhA_{k,1}, \calhA_{k,2}, \calhC_k)$ are such that $\hat A_{k,1}$ or $\hat C_k$ is independent of $\hat A_{k,2}$. Given that \ref{S2}  and \ref{S4}  do not hold for $(\calhA_{k,1}, \calhA_{k,2}, \calhC_k)$, there exist $(\hat A_{k,1}', \hat A_{k,2}') \in\calP(\calhA_{k,1},\calhA_{k,2})$ such that $\hat A_{k,1}'$ is not independent of $\hat A_{k,2}'$, and $(\hat A_{k,2}'', \hat C_{k}'') \in \calP(\calhA_{k,2},\calhC_{k})$ such that $\hat C_{k}$ is not independent of $\hat A_{k,2}''$. 
	Then we must have that  $\hat C_k$ is independent of $\hat A_{k,2}'$, and $\hat A_{k,1}'$ is independent of $\hat A_{k,2}''$. 
	We claim that there exists $\hat A_{k,2}\in\calP(\hat \calA_{k,2})$ such that both $\hat A_{k,1}'$ and $\hat C_k$ are independent of $\hat A_{k,2}$. 
	To show this, consider the matrix $\hat A_{k,2}(\alpha) = (1-\alpha)\hat A_{k,2}' + \alpha\hat A_{k,2}''$ for $\alpha\in\bbR$. Note that  $\hat A_{k,2}(\alpha)\in \calP(\calhA_{k,2})$ for all but finitely many $\alpha\in\bbR$. 
	Furthermore, $\hat A_{k,1}'$ is independent of $\hat A_{k,2}(\alpha)$ for all but finitely many $\alpha\in\bbR$. Indeed, $\hat A_{k,1}'$ is independent of $\hat A_{k,2}(\alpha)$ if and only if
	\begin{equation*}
	p(\alpha) = \det \left( \bbm \hat A_{k,1}' \\ \hat A_{k,2}(\alpha) \ebm \right)  \neq 0.
	\end{equation*}
	But $p(\alpha)$ is a polynomial in $\alpha$ and $p(1) \neq 0$, hence $p(\alpha) = 0$ only at the finitely many roots of $p$. Similarly, $\hat C_{k}$ is independent of $\hat A_{k,2}(\alpha)$ for all but finitely many $\alpha \in\bbR$ because
	\begin{equation*}
	q(\alpha) = \det \left( \bbm \hat C_k \\ \hat A_{k,2}(\alpha) \ebm \right)
	\end{equation*}
	is a polynomial in $\alpha$ and $q(0) \neq 0$. Therefore, there exists $\alpha' \in \bbR$ such that $\hat A_{k,2}(\alpha')\in\calP(\calhA_{k,2})$ and both $\hat A_{k,1}'$ and $\hat C_k$ are independent of $\hat A_{k,2}(\alpha')$. Let $\hat A_{k,2} = \hat A_{k,2}(\alpha')$ and note that $\hat A_{k,2}$ has full row rank. Since $\hat A_{k,1}'$ is independent of $\hat A_{k,2}$, it follows that
	\begin{equation*}
	\bbm \hat A_{k,1}' \\ \hat A_{k,2} \ebm
	\end{equation*}
	is nonsingular, hence there exist $y_{k,1} \in \bbR$ and $y_{k,2} \in \bbR^{\hat n_k - 1}$ such that
	\begin{equation*}
	\bbm y_{k,1} & y_{k,2}^\top \ebm \bbm \hat A_{k,1}' \\ \hat A_{k,2} \ebm = \lambda_k \hat C_k.
	\end{equation*}
	But then $y_{k,1}\neq 0$ because $\hat C_k$ is independent of $\hat A_{k,2}$ and thus $\hat A_{k,1} = \frac{y_{k,1}}{x_{k,1}} \hat A_{k,1}'\in\calP(\hat \calA_{k,1})$ is such that \eqref{eq:condition_then} holds.
	
	In conclusion, we have shown that under each of the conditions \ref{C1}-\ref{C6}, there exist $(\hat A_{k,1}, \hat A_{k,2}, \hat C_k) \in\calP(\calhA_{k,1}, \calhA_{k,2}, \calhC_k)$ and $y_{k,2} \in\bbR^{\hat n_k - 1}$ such that \eqref{eq:condition_then} holds, hence $\bbm \calA+\calB\calW\calC & \calB\calH \ebm$ has full row rank if $\bbm \calhA+\calhB\calW\calhC & \calhB\calH \ebm$ has full row rank. 
	To show that $\bbm \calA+\calB\calW\calC & \calB\calH \ebm$ has full row rank only if $\bbm \calhA+\calhB\calW\calhC & \calhB\calH \ebm$ has full row rank, just interchange the role of $(\calA_{k,1}, \calA_{k,2}, \calC_k)$ and $(\calhA_{k,1}, \calhA_{k,2}, \calhC_k)$ in the arguments above.
\end{IEEEproof}

\subsection{Proof of Lemma \ref{l:candidate_1}}
\begin{IEEEproof}
It is straightforward to verify that for $i= 1,4,6$ the triples $(\calhA_{\Ci}, \calhB_{\Ci}, \calhC_{\Ci})$ given by \eqref{eq:standard_A} and \eqref{eq:standard_BC} satisfy condition \Ci{}. Thus, for $i= 1,4,6$, the claim of the lemma follows immediately from Theorem~\ref{t:AC4equal}. Next consider condition~\Ci{} for $i=2,3,5$. 
Note that we cannot use Theorem~\ref{t:AC4equal} directly because $(\calhA_{\Ci}, \calhB_{\Ci}, \calhC_{\Ci})$ has state space dimension~1. To overcome this difficulty, introduce the auxiliary triple $(\calhA_{\Ci}',\calhB_{\Ci}',\calhC_{\Ci}')$ with
		\beqn
		\calhA_{\Ci}' = \bbm \calhA_{\Ci} & 0 \\ * & * \ebm,\quad \calhB_{\Ci}' = \bbm \calhB_{\Ci} \\ 0 \ebm,\quad \calhC_{\Ci}' = \bbm \calhC_{\Ci} & 0\ebm.
		\eeqn
		It is easily verified that $(\calhA_{\Ci}', \calhB_{\Ci}', \calhC_{\Ci}')$ satisfies condition \Ci{}. Hence, due to Theorem~\ref{t:AC4equal}, $(\calA_k, \calB_k, \calC_k)$ is equivalent to $(\calhA_{\Ci}', \calhB_{\Ci}', \calhC_{\Ci}')$. We will show that $(\calhA_{\Ci}', \calhB_{\Ci}', \calhC_{\Ci}')$ is equivalent to $(\calhA_{\Ci}, \calhB_{\Ci}, \calhC_{\Ci})$ by using Lemma~\ref{l:condition}. 
		With this in mind, suppose that $(\hat{A}_{\Ci}', \hat{B}_{\Ci}' ,\hat{C}_{\Ci}') \in \calP(\calhA_{\Ci}', \calhB_{\Ci}', \calhC_{\Ci}')$, $x\in\bbR^2$ and $\lambda\in\bbR$ are such that \[x\neq 0 \qand x^\top \hat{A}_{\Ci}' = \lambda \hat{C}_{\Ci}'.\] Then $x_2 = 0$ and $y = x_1$ is such that 
		\begin{equation*}
		y \neq 0,\quad y^\top \hat{A}_{\Ci} = \lambda \hat{C}_{\Ci},\quad y^\top\hat{B}_{\Ci}' = x^\top \hat{B}_{\Ci},
		\end{equation*}
		where $(\hat{A}_{\Ci}, \hat{B}_{\Ci},\hat{C}_{\Ci}) \in \calP(\calhA_{\Ci}, \calhB_{\Ci}, \calhC_{\Ci})$. Conversely, suppose that $(\hat{A}_{\Ci}, \hat{B}_{\Ci} ,\hat{C}_{\Ci}) \in \calP(\calhA_{\Ci}, \calhB_{\Ci}, \calhC_{\Ci})$, $y\in\bbR$ and $\lambda\in\bbR$  are such that $y \neq 0$ and $y^\top \hat{A}_{\Ci} = \lambda \hat{C}_{\Ci}$. Let
		\begin{equation*}
		\hat{A}_{\Ci}' = \bbm \hat{A}_{\Ci} & 0 \\ 1 & 1 \ebm,\quad \hat{B}_{\Ci}' = \bbm \hat{B}_{\Ci} \\ 0 \ebm,\quad \hat{C}_{\Ci}' = \bbm \hat{C}_{\Ci} & 0 \ebm.
		\end{equation*}
		Then $x = \bbm y & 0 \ebm^\top$ is such that
		\begin{equation*}
		x \neq 0,\quad x^\top \hat{A}_{\Ci}' = \lambda \hat{C}_{\Ci}',\quad x^\top\hat{B}_{\Ci}' = y^\top \hat{B}_{\Ci}.
		\end{equation*}
		From Lemma~\ref{l:condition} it follows that $(\calhA_{\Ci}', \calhB_{\Ci}', \calhC_{\Ci}')$ is equivalent to $(\calhA_{\Ci}, \calhB_{\Ci}, \calhC_{\Ci})$, hence $(\calA_k, \calB_k, \calC_k)$ is equivalent to $(\calhA_{\Ci}, \calhB_{\Ci}, \calhC_{\Ci})$ as well.
\end{IEEEproof}

\subsection{Proof of Lemma \ref{l:in_color}}
\blem \label{l:zero}
Let $\mathcal{M} \in \{0,\ast,?\}^{p \times q}$ be a pattern matrix and consider the corresponding graph $G(\calM) = (V,E_\ast\cup E_{?})$. Suppose that each node in $G(\calM)$ is colored white or black and let $D\in \mathbb{R}^{p \times p}$ be the diagonal matrix defined by
	$$
	D_{\ell \ell} = \begin{cases} 1 & \text{if node } \ell \text{ is black}, \\
	0 & \text{otherwise.}
	\end{cases}
	$$
If node $i$ has exactly one white out-neighbor $j$ and $(i,j)\in E_*$, then for all $M \in \mathcal{P}(\mathcal{M})$ and $z \in \bbR^p$ we have that $z ^\top \bbm M & D \ebm = 0$ if and only if $z^\top \bbm M & D + e_j e_j^\top \ebm = 0$, where $e_j\in\bbR^{p}$ is the $j$-th standard basis vector.
\elem

\begin{IEEEproof}
	Suppose that $M \in \calP(\calM)$ and $z \in \bbR^p$ are such that
		$z ^\top \bbm M & D  + e_j e_j^\top \ebm = 0$. 
	Clearly, $z^\top (D  + e_j e_j^\top) = 0$ implies that $z ^\top D = 0$ and thus $z^\top \begin{bmatrix} M & D \end{bmatrix} = 0$.
	
 Conversely, suppose that $M \in \calP(\calM)$ and $z \in \bbR^p$ are such that $z^\top \begin{bmatrix} M & D \end{bmatrix} = 0$. 
		Note that $z^\top \bbm M & D+ e_j e_j^\top \ebm= z^\top \bbm M & D\ebm = 0$ if $z_j = 0$, hence it is sufficient to show that $z_j = 0$. To this end, since $i$ has exactly one white out-neighbor $j$, from the $i$-th column of $z^\top M = 0$ we infer that
		\begin{equation*}
		z_jM_{ji} + \sum_{\ell \text{ is black}}z_\ell M_{\ell i} = 0.
		\end{equation*}
		But $z^\top D = 0$ implies that $z_\ell =0$ if $\ell$ is black, hence the latter reduces to $z_j M_{ji} = 0$.
		Since $(i,j)\in E_{\ast}$ we have $M_{ji} \neq 0$, which implies $z_j = 0$. This completes the proof.
\end{IEEEproof}

We are now ready to provide a Proof of Lemma \ref{l:in_color}.

\begin{IEEEproof}[Proof of Lemma \ref{l:in_color}]
	Let $\calT = \col(\calM, \calN)$ and $\calS(\calT)$ be the derived set of  $G(\calT)$.
	By applying Lemma~\ref{l:zero} repeatedly after every color change, we conclude that for all $T\in\calP(\calT)$ and $z\in\bbR^{r+1}$ it holds that $z^\top T = 0$ if and only if $z^\top \bbm T & \sum_{i \in S(\calT)} e_i e_i^\top \ebm = 0$, where $e_i \in\bbR^{r+1}$ is the $i$-th standard basis vector.

      To prove the `if' part, suppose that $1\in S(\calT)$.
	Note that $z^\top \sum_{i \in S(\calT)} e_i e_i^\top = 0$ implies that $z_i = 0$ if node $i\in S(\calT)$. 
	Therefore, for all $T\in\calP(\calT)$ and $z\in\bbR^{r+1}$ we have that $z^\top T = 0$ implies $z_1 = 0$, i.e., $\calM$ is independent of $\calN$.
	
Conversely, suppose that $1 \notin S(\calT)$. We will then show that there exist $T\in\calP(\calT)$ and $z \in\bbR^{r+1}$ such that 
\beq \label{eq:zT}
z_1 \neq 0 \qand z^\top T = 0,
\eeq
i.e., $\calM$ is not independent of $\calN$.
 Note that for all $i \in V$ and $j \notin S(\calT)$, it follows that either $(i,j) \notin E_\ast$ or $(i,j) \in E_\ast$  but there exists a node $\ell \neq j$ such that $\ell \notin S(\calT)$ and $(i,\ell) \in E_\ast\cup E_{?}$. This implies that
		\begin{equation*}
		\sum_{j \notin S(\calT)} \calT_{ji} = \begin{dcases*}
		{0} & if $\calT_{j i} = 0$ for all $j\notin S(\calT)$,\\
		{?} & otherwise.
		\end{dcases*}
		\end{equation*}
		for all $i \in V$. 
		In other words, there exists $T\in\calP(\calT)$ such that $\sum_{j\notin S(\calT)} T_{ji} = 0$ for all $i \in V$, hence the vector $z\in\bbR^{r+1}$ defined by
		\begin{equation*}
		z_j = \begin{dcases*}
		1 & if $j\not\in S(\calT)$,\\
		0 & otherwise,
		\end{dcases*}
		\end{equation*}
    		is such that \eqref{eq:zT} holds. This completes the proof.
\end{IEEEproof}

\ifCLASSOPTIONcaptionsoff
  \newpage
\fi

\bibliographystyle{IEEEtran}
\bibliography{reference}

\vfill

\end{document}